# Better Bell inequalities
# (passion at a distance)

## Richard D. Gill [1,*,†]

*Mathematical Institute, Leiden University and EURANDOM, NWO*

**Abstract:** I explain so-called quantum nonlocality experiments and discuss how to optimize them. Statistical tools from missing data maximum likelihood are crucial. New results are given on CGLMP, CH and ladder inequalities. Open problems are also discussed.

## 1. The name of the game

**QM vs. LR.** Bell's [5] theorem states that quantum physics (*aka* quantum mechanics, QM) is incompatible with classical physics. His proof exhibits a pattern of correlations, predicted in a certain situation by quantum physics, which is forbidden by any physical theory having a certain basic (and formerly uncontroversial) property called *local realism* (LR). Under LR, correlations must satisfy a *Bell inequality*, which however under QM can be violated.

Local realism = locality + realism, is closely connected to causality; a precise mathematical formulation will follow later. As we will see then, a further basic (and also uncontroversial) assumption called *freedom* needs to be made as well.

For the time being I offer the following explanatory remarks. Let us agree that the task of physics is to provide a causal explanation (or if you prefer, description) of reality. Events have causes (realism); cause and effect are constrained by time and space (locality). *Realism* has been taken for granted in physics since Aristotle; together with *locality* it has been a permanent feature and criterion of basic sanity till Einstein and others began to uncover disquieting features of quantum physics, see Einstein, Podolsky and Rosen [11], referred to hereafter as EPR.

For some, John Bell's theorem is a reason to argue that quantum physics must dramatically break down at some (laboratory accessible) level. For Bohr it would merely have confirmed the Copenhagen view that there is *no* underlying classical reality behind quantum physics, *no* Aristotelian/Cartesian/rationalist explanation of the random outcomes of quantum measurements. For others, it is a powerful incentive to deliver experimental proof that Nature herself violates local realism.

---

*This paper is dedicated to my friend Piet Groeneboom on the occasion of his 65th birthday. I started the research during my previous affiliation at the Mathematical Institute, Utrecht University. I acknowledge financial support from the European Community project RESQ, contract IST-2001-37559. The paper is based on work in progress joint with Toni Acin, Marco Barbieri, Wim van Dam, Nicolas Gisin, Peter Grünwald, Jan-Åke Larsson, Philipp Pluch, Stefan Zohren, and Marek Żukowski. Last but not least, Piet's programming assistance was vital. *Lang zal hij leven, in de gloria!*

†NWO is the Dutch national Science Foundation.

[1]Mathematical Institute, Snellius Bldg, University of Leiden, Niels Bohrweg 1, 2333 CA Leiden, Netherlands, e-mail: gill@math.leidenuniv.nl; url: http://www.math.leidenuniv.nl/~gill

By *communis opinio*, the splendid experiment of Aspect, Dalibard, and Grangier [3] settled the matter in favour of quantum physics. However, insiders have long known that that experiment has major shortcomings which imply that the matter is not settled at all. Twenty-five years later these shortcomings have *still* not been overcome, despite a continuing and intense effort and much progress; see Gill [14, 15], Santos [25]. I can report that certain experimenters think that a definitive successful experiment might well be achieved within ten years. A competition seems to be on to do it first. We will see.

**Bell-type experiments.** We are going to study the sets of all possible joint probability distributions of the outcomes of a *Bell-type experiment*, under two sets of assumptions, corresponding respectively to local realism and to quantum physics. Bell's theorem can be reformulated as saying that the set of LR probability laws is strictly contained in the QM set. But what is a Bell-type experiment?

That is not so difficult to explain. Here is a description of a $p \times q \times r$ Bell experiment, where $p$, $q$ and $r$ are fixed integers all at least equal to 2. The experiment involves a diabolical source, *Lucifer*, and a number $p$ of players or *parties*, usually called *Alice*, *Bob*, and so on. Lucifer sends a package to Alice and each of her friends by FedEx. *After* the packges have been handed over by Lucifer to FedEx, but *before* each party's package is delivered at his or her laboratory, each of the parties commits him or herself to using one particular tool or measurement-device out of some fixed set of toolboxes with which to open their packages. Suppose each party can choose one out of $q$ tools; each party's tools are labelled from 1 to $q$. There is no connection between different party's tools (and it is just for simplicity that we suppose each party has the same number). The $q$ tools of each party are conventionally called *measurements* or *settings*.

When the packages arrive, each of the parties opens their own package with the measurement setting that they have chosen. What happens precisely now is left to the reader's imagination; but we suppose that the possible outcomes for each of the parties can all be classified into one of $r$ different *outcome* categories, labelled from 0 to $r-1$. Again, there is not necessarily any connection between the outcome category labelled $x$ of different measurements for the same or different parties.

Given that Alice chose setting $a$, Bob $b$, and so on, there is some joint probability $p(x, y, \ldots | a, b, \ldots)$ that Alice will then observe outcome $x$, Bob $y$, …. We suppose that the parties chose their settings $a, b, \ldots$, at random from some joint distribution with probabilties $\pi(a, b, \ldots)$; $a, b, \ldots = 1, \ldots, q$. Altogether, one run of the whole experiment has outcome $(a, b, \ldots; x, y, \ldots)$ with probability $p(a, b, \ldots; x, y, \ldots) = \pi(a, b, \ldots) p(x, y, \ldots | a, b, \ldots)$.

If the different party's settings are independent, then each party would in practice generate their own setting in their own laboratory according to its marginal distribution. In general however we need a trusted, independent, referee, who we will call Piet, who generates the settings of all parties simultaneously and makes sure that each one receives their own setting in separate, sealed envelopes.

One can (and should) also consider "unbalanced" experiments with possibly different numbers of measurements per party, different numbers of outcomes per party's measurement. Moreover, more complicated multi-stage measurement strategies are sometimes considered. We stick here to the basic "balanced" designs, just for ease of exposition.

**The classical polytope.** Local realism and freedom can be taken mean the following:

> Measurements which were not done also have outcomes; actual and potential measurement outcomes are independent of the measurement settings actually used by all the parties.

The outcomes of measurements which were not actually done are obviously counterfactual. I am not claiming the actual existence in physical reality of these outcomes, whatever that might be supposed to mean (see EPR for one possible definition). I am supposing that a mathematical model for the experiment does allow the existence of such variables.

To argue this point, consider a computer simulation of the Bell experiment in which Lucifer's packages are put together on a classical computer, using randomization if necessary, while what goes on in each party's laboratory is also simulated on a computer. The package that is sent to each party can therefore be represented by a random number. What happens in each party's lab is the result of inputting the message from Lucifer, and the setting from Piet the referee, into another computer program which might also make use of random number generation. There can be any kind of dependence between the random numbers used in Lucifer's, Alice's, Bob's ... computers. But without loss of generality *all* this randomization might as well be done at Lucifer's computer; Alice's computer merely evaluates some function of the message from Lucifer, and the setting from Piet. We see that the outcomes are now simultaneously defined of every measurement which each party might choose, simply by considering all possible arguments to their computers programs. The assumption of freedom is simply that Piet's settings are independent of Lucifer's random numbers. Now, given Lucifer's randomization, everything that happens is completely deterministic: the outcome of each possible measurement of each party is fixed.

For ease of notation, consider briefly a two party experiment. Let $X_1, \ldots, X_q$ and $Y_1, \ldots, Y_q$ denote the counterfactual outcomes of each of Alice's and Bob's possible $q$ measurements (taking values in $\{0, \ldots, r-1\}$. We may think of these in statistical terms as missing data, in physical terms as so-called hidden variables. Denote by $A$ and $B$ Alice's and Bob's random settings, each taking values in $\{1, \ldots, q\}$. The actual outcomes observed by Alice and Bob are therefore $X = X_A$ and $Y = Y_B$. The data coming from one run of the experiment, $A, B, X, Y$, has joint probability distribution with mass function $p(a, b; x, y) = \pi(a, b, \ldots) p(x, y, |a, b) = \pi(a, b) \Pr(X_a = x, Y_b = y)$.

Now the joint probability distribution of the $X_a$ and $Y_b$ can be arbitrary, but in any case it is a mixture of all possible degenerate distributions of these variables. Consequently, for fixed setting distribution $\pi$, the joint distribution of $A, B, X, Y$ is also a mixture of the possible distributions corresponding to degenerate (deterministic) hidden variables. Since there are only finitely many degenerate distributions when $p$, $q$ and $r$ are all fixed, we see that

> Under local realism and freedom, the joint probability laws of the observable data lie in a convex polytope, whose vertices correspond to degenerate hidden variables.

We call this polytope the classical polytope.

**The quantum body.** Introductions to quantum statistics can be found in Gill [13], Barndorff-Nielsen et al. [4]. The bible of quantum information, Nielsen and Chuang [22], is a splendid resource and has introductory material for beginners to the field whether coming from physics, computer science or mathematics. The basic rule for computation of a probability distribution in quantum mechanics is called Born's law: take the squared lengths of the projections of the state vector into a collection of orthogonal subspaces corresponding to the different possible outcomes. For ease of notation, consider a two-party experiment. Take two complex Hilbert spaces $\mathcal{H}$ and $\mathcal{K}$. Take a unit vector $|\psi\rangle$ in $\mathcal{H} \otimes \mathcal{K}$. For each $a$, let $L_x^a$, $x = 0, \ldots, r-1$, denote orthogonal closed subspaces of $\mathcal{H}$, together spanning all of $\mathcal{H}$. Similarly, let $M_y^b$ denote the elements of $q$ collections of decompositions of $\mathcal{K}$ into orthogonal subspaces. Finally, define $p(x, y|a, b) = \|\Pi_{L_x^a} \otimes \Pi_{M_y^b} |\psi\rangle\|^2$, where $\Pi$ denotes orthogonal projection into a closed subspace. The reader should verify (basically by Pythagoras' theorem), that this does define a collection of joint probability distributions of $X$ and $Y$, indexed by $(a, b)$. As before we take $p(a, b, \ldots; x, y, \ldots) = \pi(a, b, \ldots) p(x, y, \ldots |a, b, \ldots)$.

The following fact is not trivial:

> The collection of all possible quantum probability laws of $A, B, X, Y$ (for fixed setting distribution $\pi$) forms a closed convex body containing the local polytope.

Beyond the $2 \times 2 \times 2$ case very little indeed is known about this convex body.

**The no-signalling polytope.**  The two convex bodies so far defined are forced to live in a lower dimensional affine subspace, by the basic normalization properties of probability distributions: $\sum_{x,y} p(a, b; x, y) = \pi(a, b)$ for all $a$, $b$. Moreover, probabilities are necessarily nonnegative, so this restricts us further to some convex polytope. However, physics (locality) implies another collection of equality constraints, putting us into a still smaller affine subspace. These constraints are called the no-signalling constraints: $\sum_y p(a, b; x, y)$ should be independent of $b$ for each $a$ and $x$, and vice versa. It is easy to check that both the local realist probability laws, and the quantum probability laws, satisfy no-signalling. Quantum mechanics is certainly a local theory as far as manifest (as opposed to hidden) variables are concerned.

> The set of probability laws satisfying no-signalling is therefore another convex polytope in a low dimensional affine subspace; it contains the quantum body, which in turn contains the classical polytope.

**Bell and Tsirelson inequalities.**  "Interesting" faces of the classical polypope, i.e., faces which do *not* correspond to the positivity constraints, generate (generalized) Bell inequalities, that is, linear combinations of the joint probabilities of the observable variables which reach a maximum value at the face. Similarly, "interesting" supporting hyperplanes to the quantum body correspond to (generalized) Tsirelson inequalities. These latter inequalities can be recast as inequalities concerning expectation values of certain observables called Bell operators.

The original Bell (more precisely, CHSH – Clauser, Horne, Shimony and Holt [6]) and Cirel'son [8] inequalities concern the $2 \times 2 \times 2$ case. However we will proceed by proving Bell's theorem – the quantum body is strictly larger than the local polytope – in the $3 \times 2 \times 2$ case for which a rather elegant proof is available due to Greenberger, Horne and Zeilinger [17].

By the way, the subtitle "passion at a distance" is a phrase coined by Abner Shimony and it expresses that though there is no action at a distance (no manifest non-locality), still quantum physics seems to allow the physical system at Alice's site to have some feeling for what is going on far away at Bob's. Rather like the oracles of antiquity, no-one can make any sense of what the oracle is saying till it is too late . . . . But one *can* use these non-classical correlations, as the physicists like to call them, to enable Alice and her friends to succeed at certain collaborative tasks, in which Lucifer is their ally while Piet is their adversary, with larger probability than is possible under any possible classical-like physics. The following example should inspire the reader to imagine such a task.

**GHZ paradox.** We consider a now famous $3 \times 2 \times 2$ example due to Greenberger, Horne and Zeillinger [17]. We use this example partly for fun, partly to exemplify the computation of Bell probability laws under quantum mechanics and under local realism.

Firstly, under local realism, one can introduce hidden variables $X_1$, $X_2$, $Y_1$, $Y_2$, $Z_1$, $Z_2$, standing for the counterfactual outcomes of Alice, Bob and Claudia's measurements when assigned settings 1 or 2 by Piet. These variables are binary, and we may as well denote their possible outcomes by $\pm 1$. Now note that

$$(X_1 Y_2 Z_2).(X_2 Y_1 Z_2).(X_2 Y_2 Z_1) = (X_1 Y_1 Z_1).$$

Thus, if the setting patterns $(1, 2, 2)$, $(2, 1, 2)$ and $(2, 2, 1)$ *always* result in $X$, $Y$ and $Z$ with $XYZ = +1$, it will also be the case the setting pattern $(1, 1, 1)$ *always* results in $X$, $Y$ and $Z$ with $XYZ = +1$.

Next define the $2 \times 2$ matrices

$$\sigma_1 \ = \ \begin{pmatrix} 0 & 1 \\ 1 & 0 \end{pmatrix}, \quad \sigma_2 \ = \ \begin{pmatrix} 1 & 0 \\ 0 & -1 \end{pmatrix}.$$

One easily checks that $\sigma_1 \sigma_2 = -\sigma_2 \sigma_1$, (anticommutation), $\sigma_1^2 = \sigma_2^2 = 1\!\!1$, the $2 \times 2$ identity matrix. Since $\sigma_1$ and $\sigma_2$ are both Hermitean, it follows that they have real eigenvalues, which by the properties given above, must be $\pm 1$.

Now define matrices $X_1 = \sigma_1 \otimes 1\!\!1 \otimes 1\!\!1$, $X_2 = \sigma_2 \otimes 1\!\!1 \otimes 1\!\!1$, $Y_1 = 1\!\!1 \otimes \sigma_1 \otimes 1\!\!1$, $Y_2 = 1\!\!1 \otimes \sigma_2 \otimes 1\!\!1$, $Z_1 = 1\!\!1 \otimes 1\!\!1 \otimes \sigma_1$, $Z_2 = 1\!\!1 \otimes 1\!\!1 \otimes \sigma_2$. It is now easy to check that

$$(X_1 Y_2 Z_2).(X_2 Y_1 Z_2).(X_2 Y_2 Z_1) = -(X_1 Y_1 Z_1),$$

and that $(X_1 Y_2 Z_2)$, $(X_2 Y_1 Z_2)$, $(X_2 Y_2 Z_1)$ and $(X_1 Y_1 Z_1)$ commute with one another.

Since these four $8 \times 8$ Hermitean matrices commute they can be simultaneously diagonalized. Some further elementary considerations lead one to conclude the existence of a simultaneous eigenvector $|\psi\rangle$ of all four, with eigenvalues $+1$, $+1$, $+1$, $-1$ respectively. We take this to be the state $|\psi\rangle$, with the three Hilbert spaces all equal to $\mathbb{C}^2$. We take the two orthogonal subspaces for the 1 and 2 measurements of Alice, Bob, and Claudia all to be the two eigenspaces of $\sigma_1$ and $\sigma_2$ respectively. This generates quantum probabilties such that the setting patterns $(1, 2, 2)$, $(2, 1, 2)$ and $(2, 2, 1)$ *always* result in $X$, $Y$ and $Z$ with $XYZ = +1$, while the setting pattern $(1, 1, 1)$ *always* results in $X$, $Y$ and $Z$ with $XYZ = -1$.

Thus we have shown that a vector of quantum probabilities exists, which cannot possibly occur under local realism. Since the classical polytope is closed, the corresponding quantum law must be strictly outside the classical polytope. It therefore violates a generalized Bell inequality corresponding to some face of the classical polytope, outside of which it must lie. It is left as an exercise to the reader to generate the corresponding "GHZ inequality."

**GHZ experiment.** This brings me to the point of the paper: how should one design good Bell experiments; and what is the connection of all this physics with mathematical statistics? Indeed there are many connections – as already alluded to, the hidden variables of a local realist theory are simply the missing data of a nonparametric missing data problem.

In the laboratory one creates the state $|\psi\rangle$, replacing Lucifer by a source of entangled photons, and the measurement devices of Alice and Bob by assemblages of polarization filters, beam splitters and photodetectors implementing hereby the measurements corresponding to the subspaces $L_a^x$, etc. One also settles on a joint setting probability $\pi$. One repeats the experiment many times, hoping to indeed observe a quantum probability law lying outside the classical polytope, i.e., violating a Bell inequality. The famous Aspect et al. [3] experiment implemented this program in the $2 \times 2 \times 2$ case, violating the so-called CHSH inequality (which we will describe later) by a large number of standard deviations. What is being done here is statistical hypothesis testing, where the null hypotheses is local realism, the alternative is quantum mechanics; the alternative being true by design of the experimenter and validity of quantum mechanics.

Dirk Bouwmeester recently carried out the GHZ experiment; the results are exciting enough to be published in *Nature* (Pan et al. [23]). He claimed in a newspaper interview that this experiment is of a rather special type: only a finite number of repetitions are necessary since the experiment exhibits events which are impossible under classical physics, but certain under quantum mechanics. However please note that the events which are certain or impossible, are only certain or impossible conditional on some other events being certain. Since the experiment is not perfect, Bouwmeester did observe some "wrong" outcome patterns, thereby destroying by his own logic the conclusion of his experiment. Fortunately his data does statistically significantly violate the accompanying GHZ inequality and publication in *Nature* was justified! The point is: *all* these experiments are statistical in nature; they do not prove for sure that local realism is false; they only give statistical evidence for this proposition; evidence which does become overwhelming if $N$, the number of repetitions, is large enough.

**How to compare different experiments.** Because of the dramatic zero-one nature of the GHZ experiment, it is felt by many physicists to be much stronger or better than experiments of the original $2 \times 2 \times 2$ CHSH type (still to be elucidated!) The original aim of the research described here was to supply objective and quantitative evaluation of such claims.

Now the geometric picture above naturally leads one to prefer an experiment where the distance from the quantum physical reality is as far as possible from the nearest local realistic or classical description. Much research has been done by physicists focussing on the corresponding Euclidean distance. However, it is not so clear what this distance means operationally, and whether it is comparable over experiments of different types. Moreover the Euclidean distance is altered by taking different setting distributions $\pi$ (though physicists usually only consider the uniform distribution). It is true that Euclidean distance is closely related to *noise resistance*, a kind of robustness to experimental imperfection. As one mixes the quantum probability distribution more and more with completely random, uniform outcomes, corresponding to pure noise in the photodetectors, the quantum probability distribution shrinks towards the center of the classical polytope, at some point passing through one of its faces. The amount of noise which can be allowed while still admitting violation of local realism is directly related to Euclidean distance, in our picture.

Van Dam, Gill and Grünwald [10] however propose to use relative entropy, $D(q : p) = \sum_{abxy} q(abxy) \log_2(q(abxy)/p(abxy))$, where $q$ now stands for the "true" probability distribution under some quantum description of reality, and $p$ stands for a local realist probability distribution. Their program is to evaluate $\sup_q \inf_p D(q : p)$ where the supremum is taken over parameters at the disposal of the experimenter (the quantum state $|\psi\rangle$, the measurement projectors, the setting distribution $\pi$; while the infimum is taken over probability distributions of outcomes given settings allowed by local realism (thus $q$ and $p$ in supremum and infimum actually stand for something different from the probability laws $q$ and $p$ lying in the quantum body and classical polytope respectively; hopefully this abuse of notation may be excused.

They argue that this relative entropy gives direct information about the number of trials of the experiment required to give a desired level of confidence in the conclusion of the experiment. Two experiements which differ by a factor 2 are such that the one with the smaller divergence needs to be repeated twice as often as the other in order to give an equally convincing rejection of local realism.

Moreover, optimizing over different sets of quantum parameters leads to various measures of "strength of non-locality." For instance, one can ask what is the best experiment based on a given entangled state $|\psi\rangle$? Experiments of different format can be compared with one another, possibly discounting the relative entropies according to the numbers of quantum systems involved in the different experiments in the obvious way (typically, a $p$ party experiment involves generation of $p$ particles at a time, so a four party experiment should be downweighted by a factor 2 when comparing with a two party experiment). We will give some examples later.

Finally, that paper showed how the interior infimum is basically the computation of a nonparametric maximum likelihood estimator in a missing data problem. Various algorithms from statistics can be succesfully applied here, in numerical rather than analytical experimentation; and progams developed by Piet Groeneboom (see Groeneboom et al. [18]) played a vital role in obtaining the results which we are now going to display.

## 2. CHSH and CGLMP

The $2 \times 2 \times 2$ case is particularly simple and well researched. In a later section, I want to compare the corresponding two particle CHSH experiment with the three particle GHZ. In another section I will discuss properties of $2 \times 2 \times d$ experiments, which form a natural generalization of CHSH and have received much attention both by theorists and experimenters in recent years. We will see that many open problems exist here and some remarkable conjectures can be posed. Preparatory to that, I will therefore now describe the so-called CGLPM inequality, the generalization from $2 \times 2 \times 2$ to $2 \times 2 \times d$ of CHSH.

For the $2 \times 2 \times d$ case an important step was made by Collins, Gisin, Linden, Massar and Popescu [9], in the discovery of a generalized Bell inequality (i.e., interesting face of the classical polytope), together with a quantum state and measurements which violated the inequality. The original specification of the inequality is rather complex, and its derivation also took two closely printed pages. Here I offer a new and extremely short derivation of an equivalent inequality, found very recently by Stefan Zohren, which further simplifies an already very simple version of my own. Proof of equivalence with the original CGLPM is tedious!

Recall that a Bell inequality is the face of a classical polytope of the form $\sum_{abxy} c_{abxy} p(abxy) \leq C$. Now since we are only concerned with probability distributions within the no-signalling polytope, the probabilities $p(abxy)$ necessarily satisfy a large number of equality constraints (normalization, no-signalling), which allows one to rewrite the Bell inequality in many different forms; sometimes remarkably different. A canonical form can be obtained by removing, by appropriate substitutions, all $p(abxy)$ with $x$ and $y$ equal to one particular value from the set of possible outcomes, e.g., outcome 0, and involving also the marginals $p(ax)$ and $p(by)$ with $x$ and $y$ non zero. This is not necessarily the "nicest" form of an inequality. However, in the canonical form the constant $C$ does disappear (becomes equal to 0).

To return to CGLMP: consider four random variables $X_1$, $X_2$, $Y_1$, $Y_2$. Note that $X_1 < Y_2$ and $Y_2 < X_2$ and $X_2 < Y_1$ implies $X_1 < Y_1$. Consequently, $X < 1 \geq Y_1$ implies $X_1 \geq Y_2$ or $Y_2 \geq X_2$ or $X_2 \geq Y_1$, and this gives us

$$\mathrm{Pr}(X_1 \geq Y_1) \;\leq\; \mathrm{Pr}(X_1 \geq Y_2) \;+\; \mathrm{Pr}(Y_2 \geq X_2) \;+\; \mathrm{Pr}(X_2 \geq Y_1).$$

This is a CGLMP inequality, when we further demand that all four variables take values in $\{0, \ldots, d-1\}$. The case $d = 2$ gives the CHSH inequality (though also in an unfamiliar form).

CGLMP describe a state and quantum measurements which generate probabilities, which violate this inequality. Take Alice and Bob's Hilbert space each to be $d$-dimensional. Consider the states $|\psi\rangle = \sum_{x=0}^{d-1} |xx\rangle/\sqrt{d}$, where $|xx\rangle = |x\rangle \otimes |x\rangle$, and $|x\rangle$ for $x = 0, \ldots, d-1$ is an orthonormal basis of $\mathbb{C}^d$. Alice and Bob's settings 1, 2 are taken to correspond to angles $\alpha_1 = 0$, $\alpha_2 = \pi/4$, and $\beta_1 = \pi/8$, $\beta_2 = -\pi/8$. When Alice or Bob receives setting $a$ or $b$, each applies the diagonal unitary operation with diagonal elements $\exp(ix\theta/d)$, $x = 0, \ldots, d-1$, to their part of the quantum system, where $\theta$ stands for their own angle (setting). Next Alice applies the quantum Fourier transform $Q$ to her part, and Bob its inverse (and adjoint) $Q^*$; $Q_{xy} = \exp(ixy/d)$, $Q^*_{xy} = \exp(-ixy/d)$. Finally Alice and Bob "measure in the computational basis", i.e., projecting onto the one-dimensional subspaces corresponding to the bases $|x\rangle$, $|y\rangle$. Applying a unitary $U$ and then measuring the projector $\Pi_M$ is of course the same as measuring the projector $\Pi_{U^*M}$; with a view to implementation in the laboratory it is very convenient to see the different measurements as actually "the same measurement" applied after different unitary transformations of each party's state have been applied. In quantum optics these operations might correspond to use of various crystals, applying an electomagnetic field across a light pulse, and so on.

That these choices gives a violation of a CGLMP inequality follows from some computation and we desperately need to understand what is going on here, as will become more obvious in a later section when I describe conjectures concerning CGLMP and these measurements.

## 3. Comparing some classical experiments: GHZ vs CHSH

First of all, let me briefly report some results from van Dam et al. [10] concerning the comparison of CHSH and GHZ. It is conjectured, and supported numerically, but not yet proved, that the best $2 \times 2 \times 2$ experiment in the sense of Kullback-Leibler divergence is the CGLMP experiment with $d = 2$ described in the last section, and usually known as the CHSH experiment. The setting probabilities should be uniform, the state is maximally entangled, the measurements are those implemented by Aspect et al. It turns out that $D$ is equal to 0.0423.... For GHZ, which is can be conjectured to be the best $3 \times 2 \times 2$ experiment, one finds $D = 0.400$, with setting probabilities uniform over the four setting patterns involved in the derivation of the paradox; zero on the other. So this experiment is apparently almost 10 times better. By the way, $D = 1$ would be the strength of the experiment when one repeatedly throws a coin which always comes up heads, in order to disprove the theory that Pr(heads) = 1/2. So GHZ is less than half as good as an experiment in which one compares probabilities 1 and 1/2; let alone comparable to an experiment comparing impossible with certain outcomes!

However in practice the GHZ experiment is not performed exactly in optimal fashion. To begin with, in order to produce each triple of photons, Bouwmeester generated two maximally entangled pairs of photons, measured the polarization of one of the four, and accepted the remaining set of three when the measured polarization was favourable, which occurs in half of the times. Since we need two pairs of photons for each triple, and discard the result half the times, the figure of merit should be divided by four. Next, the optimal setting probabilities is uniform over half of the eight possible combinations. In practice one generates settings at random at each measurement station, so that half of the combinations are actually useless. This means we have to halve again, resulting in a figure of merit for GHZ which is barely better than CHSH, and very far from the "infinity" which would correspond to an all or nothing experiment.

Actually things are even worse since the pairs of photon pairs are generated at random times and one has to be quite lucky to have two pairs generated close enough in time to one another that one has four photons to start with. Then there are the inevitable losses which further degrade the experiment . . . (more on this later). Bouwmeester needs to carry on measuring for hours in order to achieve what can be done with CHSH in minutes. Which is not to say that his experiment is not a splendid acheivement!

## 4. CGLMP as # outcomes goes to infinity

In Acin, Gill and Gisin [2] a start is made with studying optimal $2 \times 2 \times r$ experiments, and some remarkable findings were made, though almost all conclusions depend on numerics, and even on numerics depending on conjectures.

Let me first describe one rather fundamental conjecture whose truth would take us a long way in understanding what is going on.

In general nothing is known about the geometry of the classical polytope. An impossible open problem is to somehow classify all interesting faces. It is not even known if, in general, all faces which are not trivial (i.e., correspond to nonnegativity constraints) are "interesting" in the sense of being violable by quantum mechanics. As the numbers grow, the number and type of faces grow explosively, and exhausitive enumeration has only been done for very small numbers.

Clearly there are many many symmetries — the labelling of parties, measurements and outcomes is completely arbitrary. Moreover, there are three ways in which inequalities for smaller experiments remain inequalities for larger. Firstly, by merging categories in the larger experiment one obtains a smaller one, and the Bell inequalities for the smaller can be lifted to the larger. Next, by simply omitting measurements one can lift Bell inequalities for smaller experiments to larger. Finally, by conditioning on a particular outcome of a particular measurement of a particular party, one reduces a larger experiment to one with less parties, and conversely can lift a smaller inequality to a larger.

With the understanding that interesting faces for smaller polytopes can be lifted to interesting faces of larger in three different ways, the following conjecture seems highly plausible:

> All the faces of the $2 \times 2 \times r$ polytope are boring (nonnegativity) or interesting CGLMP, or lifted CGLMP, inequalities.

This is certainly true for $r = 2, 3, 4$ and 5 but beyond this there is only numerical evidence: numerical search for optimal experiments using the maximally entangled state $|\psi\rangle$ has *only* uncovered the CGLMP measurements, violating the CGLMP inequality.

Moreover this is true both using Euclidean and relative entropy distances.

The next, stunning, finding is that the best state for these experiments is not the maximally entangled state at all! Rather, it is a state of the form $\sum_x c_x |xx\rangle$ where the so-called Schmidt coefficients $c_x$ are symmetric around $x = (r-1)/2$, first decreasing and then increasing. This "U-shape" become more and more pronounced as $r$ increases. Moreover the shape is found for both figures of merit, though it is a different state for the two cases (even less entangled for divergence than for Euclidean, i.e., less entangled for statistical strength than for noise resistance). Rather thorough numerical search takes us up to about $r = 20$ and has been replicated by various researchers.

Taking as a conjecture a) that all faces are CGLMP, b) that the best measurements are also CGLMP and the state is $U$-shaped, we only need to optimize over the Schmidt coefffcients $c_x$. Numerically one can quite easily get up to about $r = 1000$ in this way. However with some tricks one can go to $r = 10\,000$ or even $100\,000$. Note that we are solving $\sup_q \inf_p D(q : p)$ where the infimum is over the local realist polytope, the supremum is just over the $c_j$. Now a solution must also be a stationary point for both optimizations. Differentiating with respect to the classical parameters, and recalling the form of $D$, one finds that one must have $\sum_{abxy} (\hat{q}_{abxy}/\hat{p}_{abxy})(p_{abxy} - \hat{p}_{abxy}) = 0$ for classical probabilities $p$ on the face of the classical polytope passing through the solution $\hat{p}$. But this face is a CGLMP inequality! Hence the coefficients, $\hat{q}_{abxy}/\hat{p}_{abxy}$ are the coefficients involved in this inequality, i.e., up to some normalization constants they are already known! However, the quantity we want to optimize, $D$ itself, is $\sum_{abxy} q_{abxy} \log_2(\hat{q}_{abxy}/\hat{p}_{abxy})$ and this is optimal over $q$ at $q = \hat{q}$ (i.e., this the accompanying Tsirelson inequality, or supporting hyperplane to the quantum body at the optimum). Since the terms in the logarithm are known (up to a normalization constant) we just have to optimize the mean of an almost known Bell operator over the state. This is a largest eigenvalue problem, numerically easy up to very very large $d$.

All this raises the question of what happens when $r \to \infty$. In particular, can one attain the largest conceivable violation of CGLMP, namely when the probability on the left is 1 and the three on the right are all 0, with infinite dimensional Hilbert spaces, and if so, are the corresponding state and measurements interesting and feasible experimentally? Strongly positive evidence and further conjectures are given in Zohren and Gill [27]. Some recent numerical results on $r = 3$ *and* 4 are given by Navascues et al. [21].

We think of this conjectured "perfect passion at a distance" as the optimal solution of a variant of the infamous game of *Polish Poker* (played in Russian bars between a Polish traveller and local Russian drinkers with the inevitable outcome that the Pole always gets the Roubles...). Now, Alice and Bob are playing together, against Piet. Piet chooses (completely randomly) a "setting" $a = 1, 2$ for Alice, and $b = 1, 2$ for Bob. Alice doesn't know Bob's setting and vice versa. Alice and Bob must now, separately, each think of a number. Denote Alice's number by $x_a$, Bob's by $y_b$. Alice and Bob's aim is to attain $x_1 < y_2$ (if Piet calls "1; 2"), and $y_2 < x_2$ (if Piet calls "2; 2"), and $x_2 < y_1$ (if ...), *and* $y_1 < x_1$ (if ...). If they choose their numbers by any classical means, e.g., with classical dice, they must fail at least a quarter of the times. However, with quantum dice (i.e., with the help of a couple of bundles of photons, donated to each of them in advance by Lucifer) they can succeed with probability arbitrarily close to certainty, by taking measurements with enough outcomes. At least, according to Zohren and Gill's conjecture...

There remains the question: *why are the CGLMP measurements optimal for the CGLMP inequality?* Where do these angles come from, what has this to do with QFT? There are some ideas about this and the problem seems ripe to be cracked.

## 5. Ladder proofs

Is the CHSH experiment the best possible experiment with two maximally entangled qubits? This seemed a very good conjecture till quite recently. However the conjecture certainly needs modification now, as I will explain.

There has been some interest recently in so-called ladder proofs of Bell's theorem. These appear to allow one to use less entangled states and get better experiments, though that dream is shown to be fallacious when one uses statistical strength as a figure of merit rather than a criterion connected to "probability zero under LR, but positive under QM" (conditional on certain other probabilities equal to zero). Exactly as for GHZ, the size of this positive probability is not very important, the experiment is about violating an inequality, not about showing that some probability is positive when it should be zero.

Let me explain the ladder idea. Consider the inequality

$$\Pr(X_1 \geq Y_1) \;\; \leq \;\; \Pr(X_1 \geq Y_2) \;\; + \;\; \Pr(Y_2 \geq X_2) \;\; + \;\; \Pr(X_2 \geq Y_1).$$

Now add to this the same inequality for another pair of hidden variables:

$$\Pr(X_2 \geq Y_2) \;\; \leq \;\; \Pr(X_2 \geq Y_3) \;\; + \;\; \Pr(Y_3 \geq X_3) \;\; + \;\; \Pr(X_3 \geq Y_2).$$

The intermediate "horizontal" 2—2 term cancels and we are left only with cross terms 1—2 and 2—3, and "end" terms 1—1 and 3—3. With a ladder built from adding four inequalities involving $X_1$ to $X_5$ and $Y_1$ to $Y_5$, out of the 25 possible comparisons, only the two end horizontal terms and eight crossing terms survive, 10 out of the total.

Numerical optimization of $D$ for longer and longer ladders, shows that actually the optimal state is always the maximally entangled state. Moreover, much to my surprise, the best $D$ is obtained with the ladder of $X_1$ to $X_5$ and $Y_1$ to $Y_5$, and it is much better than the original CHSH! However, it has a uniform distribution over 10 out of 25 combinations. If one would implement the same experiment with the uniform distribution over all 25, it becomes worse that CHSH. So the new conjecture is that CHSH is the optimal $2 \times 2 \times 2$ experiment with uncorrelated settings.

These findings come from new unpublished work with Marco Barbieri; we are thinking of actually doing this experiment.

## 6. CH for Bell

In a CHSH experiment an annoying feature is that some photons are not registered at all. This means that there are really three outcomes of each measurement, with a third outcome "no photon"; however, the outcome "no photon, no photon" is not observed at all. One has a random sample size from the conditional distribution given that there is an event in at least one of the two laboratories of Alice and Bob.

It is better to realise that the original, complete sample size is actually also random, and typically Poisson, hence the observed counts of the various events are all Poisson. But can we create useful Bell inequalities for this situation?

The answer is yes, using the possibility of reparametrization of inequalities using the equality constraints. In a $2 \times 2 \times 3$ experiment one can rewrite any Bell inequality as an inequality involving only the $p_{abxy}$ with one of $x$ or $y$ not zero, as well as the marginal probabilities $p_{ax}$, $p_{by}$ with $x$ and $y$ nonzero. The constant term in the inequality becomes 0. So one gets a linear inequality involving only observed, Poisson distributed, random variables. "Poisson statistics" allows one to supply a valid standard error even though the "total sample size" was unknown.

Applying this technique in the $2 \times 2 \times 2$ case gives a known inequality, the Clauser-Horne (CH) inequality, useful when one has binary outcomes but one of the two outcomes is not observable at all; i.e., the outcomes are "detector click" and "no detector click."

How to find a good inequality for $2 \times 2 \times 3$? I simply add a certain probability of "no event", independent on both sides of the experiment, to the quantum probabilities belonging to the classical CHSH set-up. Next I solve the problem $\inf_p D(q : p)$ using Piet Groeneboom's programs. I observe the values of $q/\hat{p}$ which define the face of the local polytope closest to $q$. I rewrite the inequality in its classical form. The result is a new inequality (not quite new: Stefano Pironio informs me it is known to N. Gisin and others) which takes account of "no event" and which is linear in the observed counts.

The linearity means that the inequality can be studied using martingale techniques to show that the experiment is "insured" against time dependence and time trends, as long as the settings are chosen randomly; cf. Gill [14, 15]. It turns out to be essentially equivalent to some rather non-linear inequalities developed by Jan-Åke Larsson, see Larsson and Gill [20], which were till now the only known way to deal with "non-events." We intend to pursue this development in the near future combining treatment of the detection, coincidence and memory loopholes (Gill [16] and Larsson and Gill [20]).

## 7. Conclusions

I did not yet mention that studying the boundary of the $2 \times 2 \times 2$ quantum body and some different generalizations led Tsirelson into some deep mathematics and connections with fundamental questions involving Grothendieck's mysterious constant, see Cirel'son [8], Tsirelson [26] (the same person . . . ), Reeds [24], and Fishburn and Reeds [12].

Bell experiments offer a rich field involving many statistical ideas, beautiful mathematics, and offering deep exciting challenges. Moreover it is a hot topic in quantum information and quantum optics. Much remains to be done.

One remains wondering why nature is like this? There are two ways nature uses to generate probabilities: one is to take a line segment of length one and cut it in two. The different experiments found by cutting it at different places are compatible with one another; one sample space will do (the unit interval). The other way of nature is to take a line segment of length one, and let it be the hypothenuse of a right angled triangle. Now the squares of the other two sides are probabilities adding to one. The different experiments are not compatible with one another (at least, in dimension three or more, according to the Kochen–Specker theorem).

According to quantum mechanics and Bell's theorem, the world is completely different from how it has been thought for two thousand years of Western science. As Vovk and Shafer recently argued, Kolmogorov was one of the first to take the radical step of associating the little omega of a probability space with the *outcome* and not the hidden *cause*. Before then, all probability in physics could be traced back to uncertainty in initial conditions. Going back far enough, one could invoke symmetry to reduce the situation to "equally likely elementary outcomes." Or more subtly, sufficient chaoticity ensures that mixed up distributions are invariant under symmetries and hence uniform. At this stage, frequentists and Bayesians use the same probabilities and get the same answers, even if they interpret their probabilities differently.

According to Bell's theorem, the randomness of quantum mechanics is truly ontological and not epistemological: it cannot be traced back to ignorance but is "for real." It is curious that the quantum physics community is currently falling under the thrall of Bayesian ideas even though their science should be telling them that the probabilities are objective. Of course, one can mix subjective uncertainties with objective quantum probabilities, but to my mind this is dissolving the baby in the bathwater, not an attractive thing to do.

Still, why is nature like this, why are the probabilities what they are? My rough feeling is as follows. Reality is discrete. Hence nature cannot be continuous. However we do observe symmetries under continuous groups (rotations, shifts); the only way to accomodate this is to make nature random, and to have the probabiltiy distributions continuous, or even covariant, with the groups. Current research in the foundations of quantum mechanics (e.g., by Inge Helland) points to the conclusions that symmetry forces the shape of the probabilities (and even forces the complex Hilbert space); just as in the Aristotelian case, but at a much deeper level, probabilities are objectively fixed by symmetries.

CGLMP inequality for infinite dimensional states. *Phys. Rev. Lett.* To appear. arxiv:quant-ph/0612020.